\documentclass[12pt]{amsart}
\usepackage{amscd,amssymb}
\newtheorem{theorem}{Theorem}

\theoremstyle{definition}
\newtheorem{definition}[theorem]{Definition}
\newtheorem{problem}[theorem]{Problem}
\newtheorem{example}[theorem]{Example}
\newtheorem{conjecture}[theorem]{Conjecture}
\theoremstyle{remark}
\newtheorem{remark}[theorem]{Remark}

\def\varph{{\varphi}}
\def\ra{{\rightarrow}}
\def\lra{{\longrightarrow}}
\def\om{{\omega}}
\def\Om{{\Omega}}

\def\al{{\alpha}}

\def\La{{\Lambda}}
\def\ga{{\gamma}}
\def\Ga{{\Gamma}}
\def\vGa{{\varGamma}}

\def\bZ{{\mathbb Z}}
\def\bQ{{\mathbb Q}}
\def\bR{{\mathbb R}}

\def\Sg{{\Sigma_g}}
\def\Mg{{\mathcal M}_g}
\def\Ig{{\mathcal I}_g}
\def\Msg{{\mathcal M}_{g,*}}
\def\Isg{{\mathcal I}_{g,*}}
\def\M1g{{\mathcal M}_{g,1}}
\def\I1g{{\mathcal I}_{g,1}}

\newcommand\Symp{\operatorname{Symp}}

\newcommand\Diff{\operatorname{Diff}}

\catcode`\@=\active 

\begin{document}

\title[Cohomological structure of the mapping class group]
{Cohomological structure of the mapping class group and beyond}
\author{Shigeyuki Morita}
\address{Department of Mathematical Sciences\\
University of Tokyo \\Komaba, Tokyo 153-8914\\
Japan}
\email{morita{\char'100}ms.u-tokyo.ac.jp}
\thanks{The author is partially supported by JSPS Grant 
16204005 and 16654011}

\begin{abstract}
In this paper, we briefly review some of the known results 
concerning the cohomological structures
of the mapping class group of surfaces, the outer automorphism group 
of free groups, the diffeomorphism group of surfaces
as well as various subgroups of them such as the Torelli group,
the $IA$ outer automorphism group of free groups, the symplectomorphism group of
surfaces. 

Based on these, we present several conjectures and
problems concerning the cohomology of these groups. 
We are particularly interested in the possible interplays between these
cohomology groups rather than merely the structures of individual groups.
It turns out that, we have to include,
in our considerations, two other groups
which contain the mapping class group as their core subgroups
and whose structures seem to be deeply related to that of the mapping class
group. They are the arithmetic mapping class group and the group of homology
cobordism classes of homology cylinders. 
\end{abstract}

\maketitle

\section{Introduction}
We begin by fixing our notations for various groups appearing
in this paper. Let $\Sigma_g$ denote a closed oriented surface
of genus $g$ which will be assumed to be greater
than or equal to $2$ unless otherwise specified.
We denote by $\mathrm{Diff}_+\Sigma_g$ the group of 
orientation preserving diffeomorphisms of $\Sg$ equipped
with the $C^\infty$ topology. The same group equipped with the
{\it discrete} topology is denoted by $\mathrm{Diff}_+^\delta\Sg$.
The mapping class group $\Mg$ is the group of path components of
$\mathrm{Diff}_+\Sigma_g$. The Torelli group, denoted by
$\Ig$, is the subgroup of $\Mg$ consisting of mapping classes
which act on the homology group $H_1(\Sg;\bZ)$ trivially.
Thus we have an extension
\begin{equation}
1\lra \Ig\lra \Mg\lra \mathrm{Sp}(2g,\bZ)\lra 1
\label{eq:tmsp}
\end{equation}
where $\mathrm{Sp}(2g,\bZ)$ denotes the Siegel modular group.
Choose an embedded disk $D\subset\Sg$ and a base point
$*\in D\subset \Sg$. We denote by $\M1g$ and $\I1g$
(resp. $\Msg$ and $\Isg$) 
the mapping class group and the Torelli group 
{\it relative} to $D$ (resp. the base point $*$).

Next let $F_n$ denote a free group of rank $n\geq 2$.
Let $\mathrm{Aut}\,F_n$ (resp. $\mathrm{Out}\,F_n$) denote
the automorphism group (resp. outer automorphism group)
of $F_n$. 
Let $\mathrm{IAut}_n$ (resp. $\mathrm{IOut}_n$) denote the subgroup of
$\mathrm{Aut}\,F_n$ (resp. $\mathrm{Out}\,F_n$)
consisting of those elements which act
on the abelianization $H_1(F_n;\bZ)$ of $F_n$ trivially.
Thus we have an extension
\begin{equation}
1\lra \mathrm{IOut}_n\lra \mathrm{Out}\,F_n\lra \mathrm{GL}(n,\bZ)\lra 1.
\label{eq:houtgl}
\end{equation}
The fundamental group $\pi_1 (\Sg\setminus \mathrm{Int}\, D)$ is a
free group of rank $2g$. Fix an isomorphism 
$\pi_1 (\Sg\setminus \mathrm{Int} D)\cong F_{2g}$.
By a classical result of Dehn and Nielsen, we can write
$$
\M1g=\{\varph\in \mathrm{Aut}\, F_{2g}, \varph(\ga)=\ga\}
$$
where the element $\ga$ is defined by
$$
\ga=[\al_1,\beta_1]\cdots [\al_g,\beta_g]
$$
in terms of appropriate free generators 
$\al_1,\beta_1,\cdots,\al_g,\beta_g$ of $F_{2g}$.
Then we have the following commutative diagram
\begin{equation}
\begin{CD}
1 @>>> \I1g @>>> \M1g @>>> 
\mathrm{Sp}(2g,\bZ) @>>> 1 \\
 @. @VVV @VVV @VVV @.\\
1 @>>> \mathrm{IAut}_{2g} @>>> \mathrm{Aut}\,F_{2g} @>>> 
\mathrm{GL}(2g,\bZ) @>>> 1.
\end{CD}
\label{eq:hautgl}
\end{equation}
Similarly, for the case of the mapping class group with respect
to a base point, we have
\begin{equation}
\begin{CD}
1 @>>> \Isg @>>> \Msg @>>> 
\mathrm{Sp}(2g,\bZ) @>>> 1 \\
 @. @VVV @VVV @VVV @.\\
1 @>>> \mathrm{IOut}_{2g} @>>> \mathrm{Out}\,F_{2g} @>>> 
\mathrm{GL}(2g,\bZ) @>>> 1.
\end{CD}
\label{eq:houtgl2}
\end{equation}

Next we fix an area form (or equivalently a symplectic form)
$\om$ on $\Sg$ and we denote by $\Symp \Sg$
the subgroup of $\Diff_+\Sg$ consisting of those elements
which preserve the form $\om$. Also let
$\Symp_0\Sg$ be the identity component of $\Symp\Sg$.
Moser's theorem \cite{Moser} implies that
the quotient group $\Symp\Sg/\Symp_0\Sg$ can be naturally
identified with the mapping class group $\Mg$ and
we have the following commutatvie diagram
\begin{equation}
\begin{CD}
1 @>>> \Symp_0\Sg @>>> \Symp\Sg @>>> 
\Mg @>>> 1 \\
 @. @VVV @VVV @| @.\\
1 @>>> \Diff_0\Sg @>>> \Diff_+\Sg @>>> 
\Mg @>>> 1
\end{CD}
\label{eq:diffsymp}
\end{equation}
where $\Diff_0\Sg$ is the identity component of $\Diff_+\Sg$.

In this paper, we also consider two other groups. Namely
the arithmetic mapping class group and the group of homology 
cobordism classes of homology cylinders. 
They will be mentioned in $\S 8$ and $\S 11$ respectively.
\par
\vspace{1cm}

\section{tautological algebra of the mapping class group}

Let $\mathbf{M}_g$ be the moduli space of smooth projective curves of genus $g$
and let $\mathcal{R}^*(\mathbf{M}_g)$ be its tautological algebra. 
Namely it is the subalgebra of the Chow algebra $\mathcal{A}^*(\mathbf{M}_g)$
generated by the tautological classes
$\kappa_i \in \mathcal{A}^{i}(\mathbf{M}_g)\ (i=1,2,\cdots)$ 
introduced by Mumford \cite{Mumford}.
Faber \cite{Faber} made a beautiful conjecture about the structure
of $\mathcal{R}^*(\mathbf{M}_g)$.
There have been done many works related to and inspired by
Faber's conjecture
(we refer to survey papers \cite{HL}\cite{Kirwan}\cite{Vakil}
for some of the recent results including enhancements of Faber's
original conjecture). However the most difficult part of Faber's 
conjecture, which claims that $\mathcal{R}^*(\mathbf{M}_g)$
should be a Poincar\'e duality algebra of dimension $2g-4$,
remains unsettled.

Here we would like to describe a topological approach to Faber's conjecture,
in particular this most difficult part. For this, we denote by
$$
e_i\in H^{2i}(\Mg;\bZ)\quad (i=1,2,\cdots)
$$
the $i$-th Mumford-Morita-Miller tautological class 
which was defined in \cite{Morita84} as follows.
For any oriented $\Sg$-bundle $\pi:E\ra X$, the 
tangent bundle along the fiber of $\pi$, denoted by $\xi$,
is an oriented plane bundle over the total space $E$. 
Hence we have its Euler class $e=\chi(\xi)\in H^2(E;\bZ)$.
If we apply the Gysin homomorphism (or the integration along the
fibers)
$\pi_*:H^{*}(E;\bZ)\ra H^{*-2}(X;\bZ)$
to the power $e^{i+1}$, we obtain a cohomology class
$$
e_i(\pi)=\pi_*(e^{i+1})\in H^{2i}(X;\bZ)
$$
of the base space $X$. By the obvious naturality of this construction,
we obtain certain cohomology classes
$$
e\in H^2(\mathrm{EDiff}_+\Sg;\bZ),\quad
e_i\in H^{2i}(\mathrm{BDiff}_+\Sg;\bZ)
$$
where $\mathrm{EDiff}_+\Sg\ra \mathrm{BDiff}_+\Sg$ denotes
the universal oriented $\Sg$-bundle.
In the cases where $g\geq 2$, 
a theorem of Earle and Eells \cite{EE} implies that the two spaces
$\mathrm{EDiff}_+\Sg$ and $\mathrm{BDiff}_+\Sg$ are
Eilenberg-MacLane spaces $K(\Msg,1)$ and $K(\Mg,1)$ 
respectively. Hence we obtain the universal
Euler class $e\in H^2(\Msg;\bZ)$ and the Mumford-Morita-Miller classes
$e_i\in H^{2i}(\Mg;\bZ)$ as group cohomology classes of the
mapping class groups.
It follows from the definition that, over the rationals, 
the class $e_i$ is the 
image of $(-1)^{i+1}\kappa_i$ 
under the natural projection
$\mathcal{A}^*(\mathbf{M}_g)\ra H^*(\Mg;\bQ)$.

Now we define $\mathcal{R}^*(\Mg)$ (resp. $\mathcal{R}^*(\Msg)$)
to be the subalgebra of $H^*(\Mg;\bQ)$ (resp. $H^*(\Msg;\bQ)$)
generated by the classes $e_1,e_2,\cdots$ 
(resp. $e,e_1,e_2,\cdots$) and call them the tautological
algebra of the mapping class group $\Mg$ (resp. $\Msg$).
There is a canonical projection 
$\mathcal{R}^*(\mathbf{M}_g)\ra\mathcal{R}^*(\Mg)$.

Let us denote simply by $H$ (resp. $H_\bQ$) 
the homology group $H_1(\Sg;\bZ)$ (resp. $H_1(\Sg;\bQ)$).
Also we set
$$
U=\La^3 H/\om_0\land H,\quad
U_\bQ=U\otimes\bQ
$$
where $\om_0\in \La^2 H$ denotes the symplectic class.
$U_\bQ$ is an irreducible representation of the algebraic group
$\mathrm{Sp}(2g,\bQ)$ corresponding to the Young diagram $[1^3]$
consisting of $3$ boxes in a single column.
Recall here that, associated to any Young diagram 
whose number of rows is less than or equal to $g$,
there corresponds an irreducible representation of $\mathrm{Sp}(2g,\bQ)$
(cf. \cite{FH}).
In our papers \cite{Morita93}\cite{Morita01}, we constructed a
morphism
\begin{equation}
\begin{CD}
\Msg @>{\rho_2}>> \Bigl(\bigl([1^2]\oplus [2^2]\bigr)
\widetilde\times_{\text{torelli}} \La^3 H_\bQ\Bigr)\rtimes \mathrm{Sp}(2g,\bQ)\\
@VVV @VVV\\
\Mg @>>{\rho_2}> \bigl([2^2]\widetilde\times U_\bQ\bigr)
\rtimes \mathrm{Sp}(2g,\bQ)
\end{CD}
\label{eq:rho2}
\end{equation}
where $[2^2]\widetilde\times U_\bQ$ denotes a central extension
of $U_\bQ$ by $[2^2]$ corresponding to the unique copy
$[2^2]\in H^2(U_\bQ)$ and 
$(\bigl([1^2]\oplus [2^2]\bigr)\widetilde\times_{\text{torelli}} \La^3 H_\bQ$
is defined similarly (see \cite{Morita01} for details).
The diagram \eqref{eq:rho2} induces the following commutative diagram.
\begin{equation}
\begin{CD}
\Bigl(\La^*\La^3 H^*_\bQ/\bigl([1^2]^{\text{torelli}}\oplus 
[2^2]\bigr)\Bigr)^{Sp} 
@>{\rho_2^*}>> H^*(\Msg;\bQ)\\
@AAA @AAA \\
\bigl(\La^* U^*_\bQ/([2^2])\bigr)^{Sp} @>>{\rho_2^*}>  H^*(\Mg;\bQ).
\end{CD}
\label{eq:sprho2}
\end{equation}
On the other hand, we proved in \cite{KM96}\cite{KM01} that the
images of the above homomorphisms $\rho_2^*$ are precisely the tautological algebras.
Here the concept of the {\it generalized} Morita-Mumford classes defined 
by Kawazumi \cite{Kawazumi98} played an important role.
Then in \cite{Morita03}, the effect of unstable
degenerations of $Sp$-modules appearing in \eqref{eq:sprho2}
was analized and in particular a part of Faber's conjecture
claiming that $\mathcal{R}^*(\Mg)$ is already generated by
the classes $e_1,e_2,\cdots,e_{[g/3]}$ was proved
(later Ionel \cite{Ionel} proved this fact at the level of 
$\mathcal{R}^*(\mathbf{M}_g)$).
Although the way of degenerations of $Sp$-modules is by no
means easy to be studied, it seems natural to expect the following.

\begin{conjecture}
The natural homomorphisms 
$$
\left(\La^*\La^3 H_\bQ\right)^{Sp}\ra H^*(\Msg;\bQ),\quad
(\La^*U_\bQ)^{Sp}\ra H^*(\Mg;\bQ)
$$ 
induce isomorphisms
\begin{align*}
\Bigl(\La^*\La^3 H^*_\bQ/\bigl([1^2]^{\text{torelli}}\oplus 
[2^2]\bigr)\Bigr)^{Sp} 
&\cong {\mathcal R}^*(\Msg)  \\
\bigl(\La^* U^*_\bQ/([2^2])\bigr)^{Sp}&\cong
{\mathcal R}^*(\Mg).
\end{align*}
Furthermore, the algebras on the left hand sides are Poincar\'e duality
algebras of dimensions $2g-2$ and $2g-4$ respectively.
\end{conjecture}

Here we mention that for a single Riemann surface $X$, the cohomology
$H^*(\mathrm{Jac}(X);\bQ)$ is a Poincar\'e duality algebra of dimension
$2g$ while it can be shown that there exists a canonical isomorphism
$$
H^*(\mathrm{Jac}(X);\bQ)/([1^2])\cong H^*(X;\bQ)
$$
which is a Poincar\'e duality algebra of dimension $2$.
Here 
$$
[1^2]\subset H^2(\mathrm{Jac}(X);\bQ)
$$
denotes
the kernel $\mathrm{Ker}(\La^2 H_\bQ^*\ra \bQ)$ of the
intersection pairing and $([1^2])$ denotes the ideal generated
by it. 
Observe that we can write
$\La^* U^*_\bQ=H^*(PH^3(\mathrm{Jac}))$ which is a 
Poincar\'e duality algebra of dimension $\binom{2g}{3}-2g$,
where $PH^3(\mathrm{Jac})$ denotes the {\it primitive part}
of the third cohomology of the Jacobian variety.
Hence the above conjecture can be rewritten as
$$
\bigl(H^*(PH^3(\mathrm{Jac}))/([2^2])\bigr)^{Sp}\cong
{\mathcal R}^*(\Mg)
$$
so that it could be 
phrased as the {\it family version} of the above simple fact
for a single Riemann surface.
\par
\vspace{1cm}

\section{Higher geometry of the mapping class group}

Madsen and Weiss \cite{MW} recently proved a remarkable result
about the homotopy type of the classifying space of
the stable mapping class group.
As a corollary, they showed that the stable rational cohomology of the 
mapping class group is isomorphic to the polynomial algebra generated
by the Mumford-Morita-Miller classes
$$
\lim_{g\to\infty} H^*(\Mg;\bQ)\cong \bQ[e_1,e_2,\cdots].
$$
We also would like to mention fundamental results of 
Tillmann \cite{Tillmann} and Madsen and Tillmann \cite{MT}.

As was explained in \cite{Morita99}, the classes $e_i$ serve as the
(orbifold) Chern classes of the tangent bundle of the moduli space
$\mathbf{M}_g$ and it may appear that, stably and quantitatively,
the moduli space $\mathbf{M}_g$ is similar to the classifying
space of the unitary group, namely the complex Grassmannian.
However, qualitatively the situation is completely different and
the moduli space has much deeper structure than the Grassmannian.
Here we would like to present a few problems concerning 
``{\it higher geometry}" of the mapping class group where we understand
the Mumford-Morita-Miller classes as the primary characteristic classes.

First we recall the following problem, because of its importance,
which was already mentioned in
\cite{Morita99} (Conjecture 3.4).

\begin{problem}
Prove (or disprove) that the even Mumford-Morita-Miller classes
$e_{2i}\in H^{4i}(\Ig;\bQ)$ are non-trivial, in a suitable stable range,
as cohomology classes of the Torelli group.
\end{problem}

The difficulty of the above problem comes from the now classical fact,
proved by Johnson \cite{Johnson80}, that the abelianization of the Torelli group
is very big, namely $H_1(\Ig;\bQ)\cong U_\bQ\ (g\geq 3)$. 
Observe that if $\Ig$ were perfect, then the above problem would have
been easily solved by simply applying the Quillen plus construction
to each group of the group extension \eqref{eq:tmsp}
and then looking at the homotopy exact sequence of the resulting
fibration.
The work of Igusa \cite{Igusa} (in particular Corollary 8.5.17)
shows a close connection between the above problem with
another very important problem (see Problem \ref{pb:Igusa} in $\S$ 4)
of non-triviality of
Igusa's higher Franz-Reidemeister torsion classes 
in $H^{4i}(\mathrm{IOut}_n;\bR)$
(Igusa uses the notation $\mathrm{Out}^h F_n$ for the group $\mathrm{IOut}_n$).
We also refer to a recent work of Sakasai \cite{Sakasai04} 
which is related to the above problem.

Next we recall the following two well-known problems about the structure of the
Torelli group which are related to a foundational work of Hain \cite{Hain97}.

\begin{problem}
Determine whether the Torelli group $\Ig\ (g\geq 3)$ is finitely
presentable or not (note that $\Ig\ (g\geq 3)$ is known to be
finitely generated by Johnson \cite{Johnson83}).
\end{problem}

\begin{problem}
Let $\mathfrak{u}_g$ denote the graded Lie algebra associated to
the prounipotent radical of the relative Malcev completion of $\Ig$ 
defined by Hain \cite{Hain97} and let 
$\mathfrak{u}_g\ra \mathfrak{h}^\bQ_g$ 
be the natural homomorphism 
(here $\mathfrak{h}^\bQ_g$ denotes the graded Lie algebra
consisting of symplectic derivations, with positive degrees, 
of the Malcev Lie algebra of $\pi_1\Sg$). 
Determine whether this homomorphism is injective or not.
\end{problem}

In \cite{Morita99}, we defined a series of secondary characteristic
classes for the mapping class group. However there was ambiguity
coming from possible odd dimensional stable cohomology classes of
the mapping class group. Because of the result of Madsen-Weiss cited
above, we can now eliminate the ambiguity and give a precise
definition as follows. For each $i$, we constructed in \cite{KM96}\cite{KM01}
explicit group cocycles $z_i \in Z^{2i}(\Mg;\bQ)$ which represent
the $i$-th Mumford-Morita-Miller class $e_i$ by making use of the
homomorphism $\Mg\ra U\rtimes \mathrm{Sp}(2g,\bZ)$ constructed
in \cite{Morita93} which extends the (first) Johnson homomorphism $\Ig\ra U$.
These cocycles are $\Mg$-invariant by the definition.
Furthermore we proved that such cocycles are unique up to coboundaries.
On the other hand, as is well known, 
any {\it odd} class $e_{2i-1}$ comes from the Siegel modular group
$\mathrm{Sp}(2g,\bZ)$ so that there is a cocycle
$z'_{2i-1}\in Z^{4i-2}(\Mg;\bQ)$ which comes from $\mathrm{Sp}(2g,\bZ)$.
This cocycle is uniquely defined up to coboundaries and $\Mg$-invariant. 
Now consider the difference $z_{2i-1}-z'_{2i-1}$. It is a coboundary so that
there exists a cochain $y_{i}\in C^{4i-3}(\Mg;\bQ)$ such that
$\delta y_i= z_{2i-1}-z'_{2i-1}$. Since $H^{4i-3}(\Mg;\bQ)=0$ by
\cite{MW} (in a suitable stable range), the cochain $y_i$
is well-defined up to coboundaries. 

Now let $\mathcal{K}_g$ be the
kernel of the Johnson homomorphism so that we have an extension
\begin{equation}
1\lra \mathcal{K}_g\lra \Ig\lra U\lra 1.
\label{eq:kg}
\end{equation}
Recall that Johnson \cite{Johnson85} proved that $\mathcal{K}_g$ 
is the subgroup of $\Mg$ generated by Dehn twists along separating
simple closed curves on $\Sg$.
The cocycle $z'_{2i-1}$ is trivial on the Torelli group $\Ig$
while the cocycle $z_{2i-1}$ (in fact any $z_i$) vanishes on $\mathcal{K}_g$.
It follows that the restriction of the cochain $y_i$ to $\mathcal{K}_g$
is a cocycle. Hence we obtain a cohomology class
$$
d_i=[y_i|_{\mathcal{K}_g}]\in H^{4i-3}(\mathcal{K}_g;\bQ).
$$
This cohomology class is $\Mg$-invariant where $\Mg$ acts on
$H^*(\mathcal{K}_g;\bQ)$ via outer conjugations. This can be shown
as follows. For any element $\varph\in\Mg$, let $\varph_*(y_i)$ be the 
cochain obtained by applying the conjugation by $\varph$ on $y_i$.
Since both cocycles $z_{2i-1}, z'_{2i-1}$ are $\Mg$-invariant,
we have $\delta \varph_*(y_i)=\delta y_i$. Hence $\varph_*(y_i)-y_i$
is a cocycle of $\Mg$. By the result of \cite{MW} again,
we see that $\varph_*(y_i)-y_i$ is a coboundary. Hence the restrictions
of $\varph_*(y_i)$ and $y_i$ to $\mathcal{K}_g$ give the same 
cohomology class.

\begin{definition}
We call the cohomology classes 
$d_i\in H^{4i-3}(\mathcal{K}_g;\bQ)^{\Mg}$
$(i=1,2 \cdots)$
obtained above the {\it secondary} characteristic classes of the mapping
class group.
\end{definition}

The secondary classes $d_i$ are stable in the following sense. Namely
the pull back of them in $H^{4i-3}(\mathcal{K}_{g,1};\bQ)$ are independent of
$g$ under natural homomorphisms induced by the inclusions 
$\mathcal{K}_{g,1}\ra\mathcal{K}_{g+1,1}$ where $\mathcal{K}_{g,1}$
denotes the subgroup of $\M1g$ generated by Dehn twists along separating
simple closed curves on $\Sg\setminus D$. This is because the
cocycles $z_{2i-1}, z'_{2i-1}$ are stable with respect to $g$.
It follows that the secondary classes $d_i$ are defined for {\it all}
$g$ as elements of $H^{4i-3}(\mathcal{K}_{g,1};\bQ)^{\M1g}$ although we have
used the result of \cite{MW}, which is valid only in a stable range.
However as elements of $H^{4i-3}(\mathcal{K}_g;\bQ)^{\Mg}$
the class $d_i$ is defined only for $g\geq 12i-9$ at present, 
although it is highly likely that it is defined for all $g$.
It was proved in \cite{Morita91} that $d_1$ is the generator
of $H^1(\mathcal{K}_g;\bZ)^{\Mg}\cong \bZ$ for all $g\geq 2$.
See \cite{Morita97} for another approach to the secondary classes.

\begin{problem}
Prove that all the secondary classes $d_2, d_3,\cdots$ are non-trivial.
\end{problem}

Here is a problem concerning the first class $d_1$. 
Let $C$ be a separating simple closed curve on $\Sg$ which
divides $\Sg$ into two compact surfaces of genera $h$ and $g-h$
and let $\tau_C\in\mathcal{K}_g$ be the Dehn twist along $C$.
Then we know that the value of $d_1$ on $\tau_C\in\mathcal{K}_g$
is $h(g-h)$ (up to a constant depending on $g$).
This is a very simple formula. However at present there is no
known algorithm to calculate the value $d_1(\varph)$ for a given
element $\varph\in\mathcal{K}_g$, say by analyzing the action of
$\varph$ on $\pi_1\Sg$. 

\begin{problem}
Find explicit way of calculating $d_1(\varph)$ for any given element
$\varph\in\mathcal{K}_g$. In particular, determine whether the
Magnus representation $\I1g\ra \mathrm{GL}(2g;\bZ[H])$ of the
Torelli group detects $d_1$ or not.
\end{problem}

Suzuki \cite{Suzuki} proved that the Magnus representation
of the Torelli group mentioned above is {\it not} faithful
so that it may happen that the intersection of the kernel of the 
Magnus representation with $\mathcal{K}_g$
is not contained in the kernel of $d_1$. We may also
ask whether the representation of the hyperelliptic mapping class group
given by Jones \cite{Jones}, restricted to the intersection of this group
with $\mathcal{K}_g$, detects $d_1$ or not (cf. Kasahara \cite{Kasahara}
for a related work for the case $g=2$). 
There are also various interesting works related to the class $d_1$
such as Endo \cite{Endo} and Morifuji \cite{Morifuji}
treating the hyperelliptic mapping class group,
Kitano \cite{Kitano} as well as Hain and Reed \cite{HR}.

Recently Biss and Farb \cite{BF} proved that the group $\mathcal{K}_g$ is
not finitely generated for all $g\geq 3$ ($\mathcal{K}_2$ is known to
be an infinitely generated free group by Mess \cite{Mess}).
However it is still not yet known whether the abelianization
$H_1(\mathcal{K}_g)$ is finitely generated or not (cf. Problem 2.2 of 
\cite{Morita99}).

Finally we would like to mention that Kawazumi \cite{Kawazumi05} is developing a 
theory of harmonic Magnus expansions which gives in particular a system of 
differential forms representing the Mumford-Morita-Miller classes
on the universal family of curves over the moduli space $\mathbf{M}_g$.

\par
\vspace{1cm}

\section{Outer automorphism group of free groups}

As already mentioned in $\S$ 1, let
$F_n$
denote a free group of rank $n\geq 2$ and let
$
{\rm Out}\, F_n = {\rm Aut}\, F_n/{\rm Inn}\, F_n
$
denote the outer automorphism group of $F_n$.
In 1986, Culler and Vogtmann \cite{CV} defined a space $X_n$,
called the {\it Outer Space}, which plays the role of the
Teichm\"uller space where the mapping class group is replaced
by ${\rm Out}\, F_n$. In particular, they proved
that $X_n$ is contractible and
${\rm Out}\, F_n$ acts on it properly discontinuously. The quotient
space
$$
\mathbf{G}_n = X_n / {\rm Out}\, F_n
$$
is called the moduli space of {\it graphs} which is the space of
all the isomorphism classes of metric graphs with fundamental group
$F_n$. Recently many works have been done on the structure of
${\rm Out}\, F_n$ as well as $\mathbf{G}_n$, notably by Vogtmann
(see her survey paper \cite{Vogtmann}),
Bestvina (see \cite{Bestvina}) and many others.

It is an interesting and important problem to compare similarity
as well as difference between the mapping class
group and ${\rm Out}\, F_n$ which will be discussed 
at several places in this book.
Here we would like to concentrate
on the cohomological side of this problem.

Hatcher and Vogtmann \cite{HV04} 
(see also Hatcher \cite{Hatcher95})
proved that the homology of 
${\rm Out}\, F_n$ stabilizes in a certain stable range.
This is an analogue of Harer's stability theorem \cite{Harer85}
for the mapping class group.
More precisely, they proved that the natural homomorphisms
$$
\mathrm{Aut}\,F_n\ra \mathrm{Aut}\,F_{n+1}, \quad
\mathrm{Aut}\,F_n\ra\mathrm{Out}\,F_n
$$ 
induce
isomorphisms on the $i$-dimensional homology group
for $n\geq 2i+2$ and $n\geq 2i+4$, respectively. 
Thus we can speak of the
stable cohomology group 
$$
\lim_{n\to\infty} \widetilde H^*(\mathrm{Out}\, F_n)
$$
of $\mathrm{Out}\,F_n$.

In the case of the mapping class group,
it was proved in \cite{Miller}\cite{Morita87}
that the natural homomorphism
${\Mg}\ra \mathrm{Sp}(2g,\bZ)$ induces an injection
$$
\lim_{g\to\infty} H^*({\rm Sp}(2g,\bZ);\bQ)
\cong \bQ[c_1,c_3,\cdots]
\subset
\lim_{g\to\infty} H^*(\Mg;\bQ)
$$
on the stable rational cohomology group
where the stable rational cohomology
of $\mathrm{Sp}(2g,\bZ)$ was determined by 
Borel \cite{Borel74}\cite{Borel81}.
In the case of ${\rm Out}\, F_n$, Igusa 
proved the following remarkable result 
which shows a
sharp contrast with the case of the mapping class group
(see Theorem 8.5.3 and Remark 8.5.4 of \cite{Igusa}).

\begin{theorem}[Igusa\cite{Igusa}]
The homomorphism
$$
\widetilde H^k({\rm GL}(n,\bZ);\bQ)\lra
\widetilde H^k({\rm Out}\, F_n;\bQ)
$$
induced by the natural homomorphism 
${\rm Out}\, F_n\ra {\rm GL}(n,\bZ)$
is the $0$-map in the stable range $n\geq 2k+1$.
\label{th:Igusa}
\end{theorem}

Recall here that the stable cohomology of 
${\rm GL}(n,\bZ)$ is given by
$$
\lim_{n\to\infty} H^*({\rm GL}(n,\bZ);\bQ)
 \cong \La_\bQ (\beta_5, \beta_9, \beta_{13},\cdots )
$$
due to Borel in the above cited papers.

On the other hand, the first non-trivial rational cohomology
of the group $\mathrm{Aut}\,F_n$ was given by Hatcher and 
Vogtmann \cite{HV}. They showed that, up to cohomology degree $6$,
the only non-trivial rational cohomology is
$$
H^4(\mathrm{Aut}\, F_4;\bQ)\cong \bQ.
$$
Around the same time, by making use of a remarkable theorem
of Kontsevich given in \cite{Kontsevich93}\cite{Kontsevich94},
the author constructed many homology classes in
$H_*(\mathrm{Out}\,F_n;\bQ)$ (see \cite{Morita99} and $\S 10$ below).
The simplest one in this construction gave a series of
elements
$$
\mu_i \in H_{4i}(\mathrm{Out}\, F_{2i+2};\bQ)
\quad (i=1,2,\cdots)
$$
and the first one $\mu_1$ was shown to be non-trivial by a
computer calculation. Responding to an inquiry of the author,
Vogtmann communicated us that she modified the argument in 
\cite{HV} to obtain an isomorphism 
$H^4(\mathrm{Out}\, F_4;\bQ)\cong \bQ$. Thus we could conclude
that $\mu_1$ is the generator of this group 
(see \cite{Morita99}\cite{Vogtmann}).
Recently Conant and Vogtmann proved that the second class 
$\mu_2\in H_8(\mathrm{Out}\, F_6;\bQ)$ is also non-trivial
in their paper \cite{CoV04} where they call $\mu_i$
the Morita classes. Furthermore they constructed many 
cycles of the moduli space $\mathbf{G}_n$ of graphs by
explicit constructions in the Outer Space $X_n$.

More recently, Ohashi \cite{Ohashi} determined the rational 
cohomology group of $\mathrm{Out}\,F_n$ for all $n\leq 6$
and in particular he showed
$$
H_8(\mathrm{Out}\, F_6;\bQ)\cong\bQ.
$$
It follows that $\mu_2$ is the generator of this group.
At present, the above two groups 
(and one more group, $H_7(\mathrm{Aut}\, F_5;\bQ)\cong\bQ$ proved by
Gerlits \cite{Gerlits}) are the only known non-trivial rational 
homology groups of $\mathrm{Out}\, F_n$ (and $\mathrm{Aut}\, F_n$).
Now we would like to present the following conjecture
based on our expectation that the classes $\mu_i$ should 
concern not only the cohomology of $\mathrm{Out}\, F_n$
but also the structure of the arithmetic mapping class group
(see $\S 8$)
as well as homology cobordism invariants of homology
$3$-spheres as will be explained in $\S 11$ below and \cite{Morita05}.

\begin{conjecture}
The classes $\mu_i$ are non-trivial for all $i=1,2,\cdots$.
\end{conjecture}

More generally we have the following.

\begin{problem}
Produce non-trivial rational (co)homology classes of $\mathrm{Out}\, F_n$.
\end{problem}

Next we consider the group $\mathrm{IOut}_n$.
In \cite{Igusa} Igusa defined higher Franz-Reidemeister torsion
classes
$$
\tau_{2i}\in H^{4i}(\mathrm{IOut}_n;\bR)
$$
as a special case of his general theory. These classes reflect
Igusa's result mentioned above (Theorem \ref{th:Igusa})
that the pull back of the Borel classes 
$\beta_{4i+1}\in H^{4i+1}(\mathrm{GL}(n,\bZ);\bR)$ 
in $H^{4i+1}(\mathrm{Out}\, F_n;\bR)$ vanish.
However it seems to be unknown whether his classes are 
non-trivial or not.

\begin{problem}[Igusa]
Prove that the higher Franz-Reidemeister torsion
classes $\tau_{2i}\in H^{4i}(\mathrm{IOut}_n;\bR)$ 
are non-trivial in a suitable stable range.
\label{pb:Igusa}
\end{problem}

In the unstable range, where the Borel classes vanish 
in $H^*(\mathrm{GL}(n,\bZ);\bR)$, there seem to be certain
relations between the classes $\tau_{2i}$, (dual of) $\mu_i$
and unstable cohomology classes in $H^*(\mathrm{GL}(n,\bZ);\bQ)$.
As the first such example, we would like to ask the following specific problem.

\begin{problem}
Prove (or disprove) that the natural homomorphism
$$
H^4(\mathrm{Out}\,F_4;\bQ)\cong\bQ \lra H^4(\mathrm{IOut}_4;\bQ)^{GL}
$$
is an isomorphism where the right hand side is generated by
(certain non-zero multiple of) $\tau_2$.
\end{problem}

Here is another very specific problem.
We know the following groups explicitly by various authors:
\begin{align*}
H^8(\mathcal{M}_{3,*};\bQ)&\cong \bQ^2\quad (\text{Looijenga \cite{Looijenga93}})\\
H^8(\mathrm{GL}(6,\bZ);\bQ)&\cong \bQ \quad 
(\text{Elbaz-Vincent, Gangl, Soul\'e \cite{EGS}})\\
H^8(\mathrm{Out}\, F_6;\bQ)&\cong \bQ \quad (\text{Ohashi \cite{Ohashi}})
\end{align*}
On the other hand, we have the following
natural injection $i$ as well as projection $p$
\begin{equation}
\mathcal{M}_{3,*}\overset{i}{\lra} \mathrm{Out}\, F_6 
\overset{p}{\lra} \mathrm{GL}(6,\bZ).
\label{eq:ip}
\end{equation}

\begin{problem}
Determine the homomorphisms 
\begin{equation}
H^8(\mathcal{M}_{3,*};\bQ) \overset{i^*}{\longleftarrow}
H^8(\mathrm{Out}\, F_6;\bQ) \overset{p^*}{\longleftarrow}
H^8(\mathrm{GL}(6,\bZ);\bQ)
\label{eq:hip}
\end{equation}
induced by the above homomorphisms in \eqref{eq:ip}.
\end{problem}

\begin{remark}
It seems to be natural to conjecture that the right map in \eqref{eq:hip}
is an isomorphism while the left map is trivial.
The former part is based on a consideration of possible geometric meaning of
the classes $\mu_i\in H_{4i}(\mathrm{Out}\, F_{2i+2};\bQ)$.
For the particular case $i=2$ here,
it was proved in \cite{EGS} that $H^9(\mathrm{GL}(6,\bZ);\bQ)=0$.
It follows that the Borel class in $H^9(\mathrm{GL}(6,\bZ);\bR)$
vanishes. Because of this, it is likely that the Igusa class
$\tau_4\in H^8(\mathrm{IOut}_6;\bR)$ would vanish as well and
the class $\mu_2$ would survive in $H_8(\mathrm{GL}(6,\bZ);\bQ)$.
For the latter part, see Remark \ref{rk:la} below.
\end{remark}

\begin{problem}
Define unstable (co)homology classes of $\mathrm{GL}(n,\bZ)$.
In particular, what can be
said about the image of $\mu_i\in H_{4i}(\mathrm{Out}\,F_{2i+2};\bQ)$
in $H_{4i}(\mathrm{GL}(2i+2,\bZ);\bQ)$
under the projection $\mathrm{Out}\,F_{2k+2}\ra \mathrm{GL}(2k+2,\bZ)$?
\end{problem}

The above known results as well as explicit computation
made so far seem to support the following conjecture
(which might be something like a folklore).

\begin{conjecture}
The stable rational cohomology of $\mathrm{Out}\, F_n$ is trivial.
Namely
$$
\lim_{n\to\infty} \widetilde H^*(\mathrm{Out}\, F_n;\bQ)=0.
$$
\end{conjecture}

We can aslo ask how the cohomology of $\mathrm{Out}\, F_n$ with 
{\it twisted coefficients} look like.

\begin{problem}
Compute the cohomology of $\mathrm{Aut}\, F_n$ and $\mathrm{Out}\, F_n$
with coefficients in various $\mathrm{GL}(n,\bQ)$-modules.
\end{problem}

For example, we could ask how Looijenga's result \cite{Looijenga96}
for the case of the mapping class group can be generalized in these
contexts. We refer to the work of Kawazumi \cite{Kawazumi05a} 
and also Satoh \cite{Sato05} for recent results 
concerning the above problem.

Finally we recall the following well known problem.

\begin{problem}
Determine whether the natural homomorphisms
\begin{align*}
&\widetilde H^*(\mathrm{Aut}\, F_{2g};\bQ)\lra \widetilde H^*(\M1g;\bQ)\\
&\widetilde H^*(\mathrm{Out}\, F_{2g};\bQ)\lra \widetilde H^*(\Msg;\bQ)
\end{align*}
induced by the inclusions 
$\M1g\ra {\rm Aut}\, F_{2g},\ \Msg\ra \mathrm{Out}\, F_{2g}$ 
are trivial or not.
\end{problem}

We refer to a result of Wahl \cite{Wahl} for a homotopy theoretical property of the
homomorphism $\M1g\ra {\rm Aut}\, F_{2g}$ where $g$ tends to $\infty$.

\begin{remark}
The known results as well as explicit computations made 
so far seem to suggest that the above maps are trivial.
According to a theorem of Kontsevich \cite{Kontsevich93}\cite{Kontsevich94}
(see $\S 9$  below), the triviality of the second map above is equivalent to 
the statement that the natural inclusion
$$
\mathfrak{l}^+_\infty \lra \mathfrak{a}^+_\infty
$$
between two infinite dimentional Lie algebras (see $\S 9$ for 
the definition) induces the {\it trivial map}
$$
H^*(\mathfrak{a}^+_\infty)^{Sp}\lra H^*(\mathfrak{l}^+_\infty)^{Sp}
$$
in the $Sp$-invariant cohomology groups.
Here the trivial map means that it is the $0$-map
except for the bigraded parts which correspond to the $0$-dimensional
homology groups of  
$\Msg$ and $\mathrm{Out}\, F_{2g}$.
\label{rk:la}
\end{remark}

\begin{remark}
In this paper, we are mainly concerned with the rational cohomology group of the 
mapping class group, $\mathrm{Out}\,F_n$ and other groups.
As for cohomology group with finite coefficients or
torsion classes, here we only mention
the work of Galatius \cite{Galatius} which determines the mod $p$
stable cohomology of the mapping class group and also
Hatcher's result \cite{Hatcher95} that the stable homology of 
$\mathrm{Out}\,F_n$ contains the homology of $\Om^{\infty}S^\infty$
as a direct summand.
\end{remark}

\par
\vspace{1cm}

\section{The derivation algebra of free Lie algebras and the traces}

As in $\S 1$, let $F_n$ be a free group of rank $n\geq 2$
and let us denote the abelianization $H_1(F_n)$ of $F_n$
simply by $H_n$. Also let $H_n^\bQ=H_n\otimes\bQ$.
Sometimes we omit $n$ and we simply write $H$ and $H_\bQ$ instead of
$H_n$ and $H_n^\bQ$.
Let 
$$
\mathcal{L}_n=\oplus_{k=1}^\infty \mathcal{L}_n(k)
$$
be the free graded Lie algebra generated by $H_n$. Also let
$\mathcal{L}_n^\bQ=\mathcal{L}_n\otimes\bQ$. 
We set
$$
\mathrm{Der}^+(\mathcal{L}_n)=\{\text{derivation}\ D
\ \text{of}\ \mathcal{L}_n\ \text{with positive degree}\}
$$
which has a natural structure of a graded Lie algebra over $\bZ$. 
The degree $k$ part of this graded Lie algebra can be expressed as
$$
\mathrm{Der}^+(\mathcal{L}_n)(k)=\mathrm{Hom}(H_n,\mathcal{L}_n(k+1))
$$
and we have
$$
\mathrm{Der}^+(\mathcal{L}_n)=\bigoplus_{k=1}^\infty 
\mathrm{Der}^+(\mathcal{L}_n)(k).
$$
Similarly we consider
$\mathrm{Der}^+(\mathcal{L}_n^\bQ)$ which is a graded Lie algebra over $\bQ$.

In the case where we are given an identification 
$\pi_1(\Sg\setminus\mathrm{Int} D)\cong F_{2g}$,
we have the symplectic class $\om_0\in \mathcal{L}_{2g}(2)=\La^2 H_{2g}$
and we can consider the following graded Lie subalgebra
\begin{align*}
\mathfrak{h}_{g,1}&=\{D\in \mathrm{Der}^+(\mathcal{L}_{2g}); D(\om_0)=0\}\\
&=\bigoplus_{k=1}^\infty \mathfrak{h}_{g,1}(k).
\end{align*}
Similarly we have $\mathfrak{h}_{g,1}^\bQ=\mathfrak{h}_{g,1}\otimes\bQ$.

In our paper \cite{Morita93a}, for each $k$, we defined a certain homomorphism
$$
\mathrm{trace}(k):\mathrm{Der}^+(\mathcal{L}_n)(k)\lra S^k H_n
$$ 
where $S^k H_n$ denotes the $k$-th symmetric power of $H_n$. We call this
``trace" because it is defined as the usual trace of 
the abelianized
{\it non-commutative Jacobian matrix} of each homogeneous derivation.
Here we recall the definition briefly from the above cited paper.
Choose a basis $x_1,\cdots,x_n$ of $H_n=H_1(F_n;\bZ)$.
We can consider $\mathcal{L}_n(k+1)$ as a natural submodule of
$H_n^{\otimes (k+1)}$ consisting of all the Lie polynomials of degree $k+1$.
For example $[x_1,x_2]\in \mathcal{L}_n(2)$ corresponds to
the element $x_1\otimes x_2-x_2\otimes x_1\in H^{\otimes 2}$.
By using the concept of the Fox free differential, we can also embed
$\mathcal{L}_n(k+1)$ into the set $(H_n^{\otimes k})^{n}$ of
all the $n$-dimensional column vectors with entries in $H_n^{\otimes k}$
by the following correspondence
$$
\mathcal{L}_n(k+1)\ni \eta\longmapsto 
 \left(\frac{\partial \eta}{\partial x_i}\right) \in (H_n^{\otimes k})^{n}.
$$
Here for each monomial $\eta\in \mathcal{L}_n(k+1)\subset H_n^{\otimes (k+1)}$ 
which is uniquely
expressed as
$$
\eta=\eta_1\otimes x_1+\cdots + \eta_n\otimes x_n\quad (\eta_i\in H_n^{\otimes k}),
$$
we have
$$
\frac{\partial \eta}{\partial x_i}=\eta_i.
$$

\begin{definition}
In the above terminologies, the $k$-th trace 
$\mathrm{trace}(k):\mathrm{Der}^+(\mathcal{L}_n)(k)\lra S^k H_n$
is defined by
$$
\mathrm{trace}(k)(f)=\left(\sum_{i=1}^n \frac{\partial f(x_i)}{\partial x_i}\right)
^{\mathrm{ab}}
$$
where $f\in \mathrm{Der}^+(\mathcal{L}_n)(k)=\mathrm{Hom}(H_n,\mathcal{L}_n(k+1))$ 
and the superscript
$\mathrm{ab}$ denotes the natural projection $H_n^{\otimes k}\ra S^k H_n$.
\end{definition}

\begin{remark}
If we identify the target $\mathrm{Hom}(H_n,\mathcal{L}_n(k+1))$ of
$\mathrm{trace}(k)$ with 
$$
H_n^*\otimes \mathcal{L}_n(k+1)\subset H_n^*\otimes H_n^{\otimes (k+1)}
$$
where $H_n^*=\mathrm{Hom}(H_n,\bZ)$ denotes the dual space of $H_n$,
then it follows immediately
from the definition that $\mathrm{trace}(k)$ is equal to the restriction of the
contraction
$$
C_{k+1}:H_n^*\otimes H_n^{\otimes (k+1)}\lra H_n^{\otimes k}
$$
followed by the abelianization $H_n^{\otimes k}\ra S^k H_n$.
Here
$$
C_{k+1}(f\otimes u_1\otimes\cdots\otimes u_{k+1})=f(u_{k+1}) u_1\otimes\cdots
\otimes u_k
$$
for $f\in H_n^*, u_i\in H_n$. Also it is easy to see that, if we replace
$C_{k+1}$ with $C_{1}$ defined by
$$
C_{1}(f\otimes u_1\otimes\cdots\otimes u_{k+1})=f(u_{1}) u_2\otimes\cdots
\otimes u_{k+1}
$$
in the above discussion,
then we obtain $(-1)^k \mathrm{trace}(k)$.
\end{remark}

\begin{example}
Let $\mathrm{ad}_{x_2}(x_1)^k\in \mathrm{Der}^+(\mathcal{L}_n)(k)$ be the element
defined by
\begin{align*}
\mathrm{ad}_{x_2}(x_1)^k(x_2)&=[x_1,[x_1,[\cdots,[x_1,x_2]\cdots]
\quad (\text{$k$-times $x_1$})\\
\mathrm{ad}_{x_2}(x_1)^k(x_i)&=0\quad (i\not=2).
\end{align*}
Then a direct computation shows that
$$
\mathrm{trace}(k)(\mathrm{ad}_{x_2}(x_1)^k)=x_1^k.
$$
\end{example}

As was mentioned in \cite{Morita93a}, the traces are $GL(H_n)$-equivariant
in an obvious sense.
Since $x_1^k$ generates $S^k H_n$
as a $GL(H_n)$-module, the above example implies that
the mapping 
$\mathrm{trace}(k):\mathrm{Der}^+(\mathcal{L}_n)(k)\lra S^k H_n$
is surjective for any $k$.
Another very important property of the traces
proved in the above cited paper is that they vanish identically
on the commutator ideal 
$[\mathrm{Der}^+(\mathcal{L}_n),\mathrm{Der}^+(\mathcal{L}_n)]$.
Hence we have the following {\it surjective}
homomorphism of graded Lie algebras
$$
(\tau_1,\oplus_k \mathrm{trace}(k)):\mathrm{Der}^+(\mathcal{L}_n)\lra 
\mathrm{Hom}(H_n,\La^2 H_n)
\oplus \bigoplus_{k=2}^\infty S^k H_n
$$
where the target is understood to be an {\it abelian} Lie algebra.
We have also proved that, for any $k$, $\mathrm{trace}(2k)$ vanishes identically
on $\mathfrak{h}_{g,1}$ and that 
$\mathrm{trace}(2k+1):\mathfrak{h}_{g,1}^\bQ(2k+1)\ra S^{2k+1}H_{2g}^\bQ$
is surjective. Thus we have a surjective homomorphism
\begin{equation}
(\tau_1,\oplus_k \mathrm{trace}(2k+1)):\mathfrak{h}_{g,1}^\bQ\lra 
\La^3 H_{2g}^\bQ
\oplus \bigoplus_{k=1}^\infty S^{2k+1} H_{2g}^\bQ
\label{eq:ttrace}
\end{equation}
of graded Lie algebras which we conjectured to give the
{\it abelianization} of the Lie algebra
$\mathfrak{h}_{g,1}^\bQ$ (see Conjecture 6.10 of \cite{Morita99}).

Recently Kassabov \cite{Kassabov} (Theorem 1.4.11) proved the following 
remarkable result.
Let $x_1,\cdots,x_n$ be a basis of $H_n$ as before.

\begin{theorem}[Kassabov]
Up to degree $n(n-1)$,
the graded Lie algebra $\mathrm{Der}^+(\mathcal{L}_n^\bQ)$ is generated as 
a Lie algebra and $\mathfrak{sl}(n,\bQ)$-module by the elements
$\mathrm{ad}(x_1)^k \ (k=1,2,\cdots)$ and the element
$D$ which sends $x_1$ to $[x_2,x_3]$ and $x_i (i\not=1)$
to $0$.
\label{th:Kassabov}
\end{theorem}

If we combine this theorem with the concept of the traces, we obtain
the following.

\begin{theorem}
The surjective Lie algebra homomorphism
$$
(\tau_1,\oplus_k \mathrm{trace}(k)):
\mathrm{Der}^+(\mathcal{L}^\bQ_n) \lra \mathrm{Hom}(H_n^\bQ,\La^2 H_n^\bQ)
\oplus \bigoplus_{k=2}^\infty S^k H_n^\bQ,
$$
induced by the degree $1$ part and the traces,
gives the abelianization of the graded Lie algebra
$\mathrm{Der}^+(\mathcal{L}^\bQ_n)$ up to degree $n(n-1)$
so that any element of degree $2\leq d\leq n(n-1)$
with vanishing trace belongs to the commutator ideal
$[\mathrm{Der}^+(\mathcal{L}^\bQ_n), \mathrm{Der}^+(\mathcal{L}^\bQ_n)]$.
Furthermore, any $\mathfrak{sl}(n,\bQ)$-equivariant
splitting to this abelianization generates
$\mathrm{Der}^+(\mathcal{L}^\bQ_n)$ in this range.
Hence stably there exists an isomorphism
$$
H_1\left(\mathrm{Der}^+(\mathcal{L}^\bQ_\infty)\right)\cong 
\mathrm{Hom}(H_\infty^\bQ,\La^2 H_\infty^\bQ)
\oplus \bigoplus_{k=2}^\infty S^k H_\infty^\bQ
$$
and the degree $1$ part and (any $\mathfrak{sl}(n,\bQ)$-equivariant
splittings of) the traces generate 
$\mathrm{Der}^+(\mathcal{L}^\bQ_\infty)$.
\end{theorem}

Although the structure of $\mathfrak{h}_{g,1}^\bQ$ is much more complicated
than that of $\mathrm{Der}^+(\mathcal{L}^\bQ_n)$,
Kassabov's argument adapted to this case together with some additional idea  
will produce enough information about the generation as
well as the abelianization of $\mathfrak{h}_{g,1}^\bQ$
in a certain {\it stable range}.
Details will be given in our forthcoming paper \cite{Morita05}.
It follows that any element in the Lie algebra $\mathfrak{h}_{\infty,1}^\bQ$
can be expressed in terms of the degree $1$ part and the traces.

\par
\vspace{1cm}

\section{The second cohomology of $\mathfrak{h}_{g,1}^\bQ$}

In this section, we define a series of elements in
$H^2(\mathfrak{h}_{g,1}^\bQ)^{Sp}$ which denotes the 
{\it $Sp$-invariant part} of the second cohomology of the 
graded Lie algebra $\mathfrak{h}_{g,1}^\bQ$.

As is well known, $U_\bQ=\La^3 H_\bQ/\om_0\land H_\bQ$ and
$S^{2k+1} H_\bQ\ (k=1,2,\cdots)$ are all irreducible representations
of $\mathrm{Sp}(2g,\bQ)$. It is well known in the representation
theory that, for any irreducible representation $V$
of the algebraic group $\mathrm{Sp}(2g,\bQ)$,
the tensor product $V\otimes V$ contains a unique trivial 
summand $\bQ\subset V\otimes V$
(cf. \cite{FH} for generalities
of the representations of the algebraic group $\mathrm{Sp}(2g,\bQ)$). 
In our case where $V$ is any of the above irreducible representations,
it is easy to see that the trivial summand
appears in the second exterior power part $\La^2 V\subset V\otimes V$.
It follows that each of
$$
\La^2 U_\bQ,\ \La^2 S^3 H_\bQ,\ \La^2 S^5 H_\bQ,\ \cdots
$$
contains a unique trivial summand $\bQ$.
Let 
$$
\iota_1: \La^2 U_\bQ\ra \bQ,\quad \iota_{2k+1}:\La^2 S^{2k+1} H_\bQ
\ra \bQ\quad (k=1,2,\cdots)
$$
be the unique (up to scalars) $Sp$-equivariant homomorphism.
We would like to call them {\it higher intersection pairing}
on surfaces which generalize the usual one
$\La^2 H_\bQ \ra \bQ$.
We can write
$$
\iota_1\in H^2(U_\bQ)^{Sp}, \quad \iota_{2k+1}\in H^2(S^{2k+1}H_\bQ)^{Sp}.
$$

\begin{definition}
We define the cohomology classes 
$$
e_1, t_3, t_5,\cdots \in H^2(\mathfrak{h}_{g,1}^\bQ)^{Sp}
$$
by setting
$$
e_1= \bar\tau_1^*(\iota_1),\quad t_{2k+1}=\mathrm{trace}(2k+1)^*(\iota_{2k+1})
$$
where $\bar\tau_1$ denotes the composition $\mathfrak{h}^\bQ_{g,1} \ra \La^3 H_\bQ
\ra U_\bQ$.
\label{def:et}
\end{definition}

\begin{conjecture}
The classes $e_1, t_3, t_5,\cdots$ are all non-trivial. Furthermore they are
linearly independent and form a basis of $H^2(\mathfrak{h}_{g,1}^\bQ)^{Sp}$.
\label{conj:t}
\end{conjecture}

\begin{remark}
The element $e_1$ is the {\it Lie algebra version} of the first
Mumford-Morita-Miller class (we use the same notation).
\end{remark}

\begin{remark}
The Lie algebra $\mathfrak{h}^\bQ_{g,1}$ is graded so that the cohomology
group $H^2(\mathfrak{h}_{g,1}^\bQ)^{Sp}$ is {\it bigraded}.
Let $H^2(\mathfrak{h}_{g,1}^\bQ)^{Sp}_{n}$ denote the weight $n$ part 
of $H^2(\mathfrak{h}_{g,1}^\bQ)^{Sp}$
(see $\S 9$ for more details). Then by definition we have
$$
e_1\in H^2(\mathfrak{h}_{g,1}^\bQ)^{Sp}_{2},
\quad t_{2k+1}\in H^2(\mathfrak{h}_{g,1}^\bQ)^{Sp}_{4k+2}.
$$
Hence if the elements $e_1, t_3, t_5,\cdots$ are non-trivial, then they are
automatically linearly independent. Thus the above conjecture can be  
rewritten as
$$
H^2(\mathfrak{h}_{g,1}^\bQ)^{Sp}_{n}\cong
\begin{cases}
\bQ\quad &(n=2,6,10,14,\cdots)\\
0 \quad &(\text{otherwise})
\end{cases}
$$
where the summands $\bQ$ in degrees $2,6,10,\cdots$ are generated by 
the above classes.
\end{remark}

As for the non-triviality, all we know at present is the non-triviality
of $e_1, t_3, t_5$. The non-triviality of the class $t_{2k+1}$
is the same as that of the class $\mu_k$
because of the theorem of Kontsevich described in $\S 9$.
This will be explained in that section.

\par
\vspace{1cm}

\section{Constructing  cohomology classes of $\mathfrak{h}_{g,1}^\bQ$}

In this section, we describe a general method of constructing 
$Sp$-invariant cohomology classes of the Lie algebra $\mathfrak{h}_{g,1}^\bQ$
which generalize the procedure given in the previous section.
As was already mentioned in our paper \cite{Morita99},
the homomorphism \eqref{eq:ttrace} induces the following homomorphim
in the $Sp$-invariant part of the cohomology
\begin{equation}
H^*\Big(\La^3 H_\bQ
\oplus \bigoplus_{k=1}^\infty S^{2k+1} H_{2g}^\bQ\Big)^{Sp}
\lra
H^*(\mathfrak{h}_{g,1}^\bQ)^{Sp}.
\label{eq:chtrace}
\end{equation}

By the same way as in \cite{Morita96}\cite{KM96},
the left hand side can be computed by certain
polynomial algebra 
$$
\bQ[\Ga; \Ga\in \mathcal{G}^{odd}]
$$
generated by graphs belonging to $\mathcal{G}^{odd}$
which denotes the set of all isomorphism classes of 
{\it connected graphs} with valencies in the set 
$$
3,3,5,7,\cdots
$$ 
of odd integers. Here we write two copies of
$3$ because of different roles: the first one is related to the 
target $\La^3 H_\bQ$ of $\tau_1$ (alternating) while the second one is related
to the target $S^3 H_\bQ$ of $\mathrm{trace}(3)$ (symmetric).
The other $2k+1\ (k=2,3,\cdots)$ are related to the target $S^{2k+1}H_\bQ$
of $\text{trace}(2k+1)$.
Thus we obtain a homomorphism
\begin{equation}
\Phi: \bQ[\Ga; \Ga\in \mathcal{G}^{odd}]\lra H^*(\mathfrak{h}_{g,1}^\bQ)^{Sp}.
\label{eq:cgh}
\end{equation}
The elements $e_1, t_3, t_5, \cdots$ defined in Definition \ref{def:et}
arise as the images, under $\Phi$, of those graphs
with exactly two vertices which are connected by
$3, 3, 5, 7, \cdots$ edges.

\begin{remark}
As was mentioned already in the previous section
$\S 6$, the cohomology of $\mathfrak{h}_{g,1}^\bQ$
is bigraded. Let $\Ga\in \mathcal{G}^{odd}$ be a connected graph
whose valencies are $v^a_3$ times $3$ (alternating), $v^s_3$ times $3$
(symmetric) and $v_{2k+1}$ times $2k+1\ (2k+1>3)$. Then
\begin{equation}
\Phi(\Ga)\in H^{d}(\mathfrak{h}_{g,1}^\bQ)^{Sp}_{n}
\label{eq:phiga}
\end{equation}
where
\begin{align*}
d&=v^a_3+v^s_3+v_5+v_7+\cdots\\
n&=v^a_3+3 v^s_3+5 v_5+7 v_7+\cdots.
\end{align*}
Observe that $n+2 v_3^a$ is equal to twice of the number of
edges of $\Gamma$. It follows that $n$ (and hence $d$) is always
an even integer.
\end{remark}

\begin{problem}
Find explicit graphs $\Ga\in \mathcal{G}^{odd}$ such that
the corresponding homology classes $\Phi(\Ga)$ are non-trivial.
\label{pb:grphi}
\end{problem}
\par
\vspace{1cm}

\section{Three groups beyond the mapping class group}

In view of the definition of the Lie algebra $\mathfrak{h}_{g,1}$ (see $\S 5$),
we may say that it is the ``Lie algebra version" of the 
mapping class group $\M1g$. However the result of the author \cite{Morita93a}
showed that it is too big to be considered so and the
following question arose: what is the algebraic and/or 
geometric meaning of the complement of the image of $\M1g$
in $\mathfrak{h}_{g,1}$ ?
Two groups came into play in this framework in the 1990's.
One is the arithmetic mapping class group through the works
of number theorists, notably Oda, Nakamura and Matsumoto,
and the other is $\mathrm{Out}\,F_n$
via the theorem of Kontsevich described in the next section $\S 9$.
In this section, we would like to consider the former group briefly
from a very limited point of view
(see \cite{Nakamura97}\cite{Matsumoto00} and references therein for details).
The latter group was already introduced in $\S 4$ and will be
further discussed in $\S 10$.

More recently, it turned out that we have to treat one more group in the above 
setting and that is the group of homology cobordism classes of
homology cylinders. This will be discussed in $\S 11$ below.
We strongly expect that the traces will give rise to meaningful invariants
in each of these three groups beyond the mapping classs group.

Now we consider the first group above.
The action of $\M1g$ on the lower central series of 
$\vGa=\pi_1 (\Sg\setminus \mathrm{Int}\, D)$ induces a filtration
$\{\M1g(k)\}_k$ on $\M1g$ as follows.
Let $\vGa_1=[\vGa,\vGa]$ be the commutator subgroup of $\vGa$
and inductively define $\vGa_{k+1}=[\vGa,\vGa_{k}]\ (k=1,2,\cdots)$.
The quotient group $N_k=\vGa/\vGa_k$ is called the $k$-th nilpotent
quotient of $\vGa$. Note that $N_1$ is canonically isomorphic to 
$H_{2g}=H_1(\Sg;\bZ)$.
Now we set
$$
\M1g(k)=\{\varph\in\M1g; \text{$\varph$ acts on $N_k$ trivially}\}.
$$
Thus the first group $\M1g(1)$ in this filtration is nothing but the
Torelli group $\I1g$. As is well known, the quotient group
$\vGa_k/\vGa_{k+1}$ can be identified with the $(k+1)$-st term
$\mathcal{L}_{2g}(k+1)$ of the free graded Lie algebra generated
by $H_{2g}$ (see $\S$ 5). It can be checked that the correspondence
\begin{align*}
\M1g(k)&\ni \varph\longmapsto \\
& \vGa\ni\alpha\mapsto \varph_*(\alpha)\alpha^{-1}\in \vGa_k/\vGa_{k+1}\cong
\mathcal{L}_{2g}(k+1)
\end{align*}
descends to a homomorphism
$$
\tau_k: \M1g(k)\lra \mathrm{Hom}(H_{2g},\mathcal{L}_{2g}(k+1))
$$
which is now called the $k$-th Johnson homomorphism
because it was introduced by Johnson (see \cite{Johnson80}\cite{Johnson83a}).
Furthermore it turns out that 
the totality $\{\tau_k\}_k$ of these homomorphims
induces an injective homomorphism
of graded Lie algebras
\begin{equation}
\mathrm{Gr}^+(\M1g)=\bigoplus_{k=1}^\infty
\M1g(k)/\M1g(k+1)\lra \mathfrak{h}_{g,1}
\label{eq:johnson}
\end{equation}
(see \cite{Morita93a}\cite{Morita99}
for details).
Although there have been obtained many important results concerning
the image of the above homomorphism, the following is still open.

\begin{problem}
Determine the image as well as the cokernel of the homomorphism \eqref{eq:johnson}
explicitly. 
\end{problem}

Note that Hain \cite{Hain97} proved that the image of \eqref{eq:johnson},
after tensored with $\bQ$, is precisely the Lie subalgebra generated by
the degree $1$ part. However it is unclear which part of 
$\mathfrak{h}_{g,1}^\bQ$ belongs to this Lie subalgebra.

In relation to this problem,
Oda predicted, in the late 1980's, that there should arise 
``arithmetic obstructions" to the surjectivity of Johnson
homomorphism.  More precisely, based on the theory of Ihara in number theory
which treated mainly the case $g=0, n=3$, he expected that
the absolute Galois group
$\mathrm{Gal}(\overline\bQ/\bQ)$ should ``appear" in 
$\mathfrak{h}_{g,1}\otimes \bZ_\ell$ outside of the
geometric part and which should be $Sp$-invariant 
for any genus $g$ and for any prime $\ell$.
In 1994, Nakamura \cite{Nakamura96} proved, among other
results, that this is in fact the case 
(see also Matsumoto \cite{Matsumoto96}).
This was the second
obstruction to the surjectivity of Johnson homomorphism,
the first one being the traces in \cite{Morita93a}.
This raised the following problem.

\begin{problem}
Describe the Galois images in $\mathfrak{h}_{g,1}\otimes\bZ_\ell$.
\end{problem}
 
The above result was proved by analyzing the number theoretical
enhancement of the Johnson homomorphism where the geometric
mapping class group is replaced by
the arithmetic mapping class group which is expressed as
an extension
$$
1\lra \hat {\mathcal M}_g^n\lra \pi_1^{alg} {\bf M}_g^n/\bQ
\lra {\rm Gal}(\overline\bQ/\bQ)\lra 1
$$
and studied by Grothendieck, Deligne, Ihara, Drinfel'd and many number
theorists.
Nakamura continued to study 
the structure of the mapping class group from the point of view of
number theory extensively (see e.g. \cite{Nakamura99}\cite{Nakamura02}).
On the other hand, Hain and Matsumoto recently
proved remarkable results concerning this subject
(see \cite{HM03}\cite{HM05}).
In view of deep theories in number theory, 
as well as the above explicit results,
it seems to be conjectured that there should exist
an embedding 
$$
\text{FreeLie}_\bZ (\sigma_3,\sigma_5,\cdots )\ \subset \mathfrak{h}_{g,1}
$$
of certain free graded Lie algebra over $\bZ$
generated by certain elements $\sigma_3,\sigma_5,\cdots$, 
corresponding to the Soul\'e elements,
into $\mathfrak{h}_{g,1}$ such that the tensor product of it with $\bZ_\ell$
coincides with the image of $\mathrm{Gal}(\overline\bQ/\bQ)$
for any prime $\ell$.

We expect that the above conjectured
free graded Lie algebra (the motivic Lie algebra) can be realized
inside $\mathfrak{h}_{g,1}$
(in fact inside the commutator ideal
$[\mathfrak{h}_{g,1}, \mathfrak{h}_{g,1}]$) 
explicitly in terms of the traces.
In some sense, the elements $\sigma_{2k+1}$
should be {\it decomposable} in higher genera.
Here we omit the precise form of the expected
formula which will be given in a forthcoming paper. 

Finally we mention the analogue of the Johnson homomorphisms for the
group $\mathrm{Aut}\,F_n$ very briefly. Prior to the 
work of Johnson, Andreadakis \cite{Andreadakis} introduced and studied
the filtration on $\{\mathrm{Aut}\,F_n(k)\}_k$ which is induced from the action
of $\mathrm{Aut}\,F_n$ on the lower central series of $F_n$.
The first group $\mathrm{Aut}\,F_n(1)$ in this filtration is nothing but the
group $\mathrm{IAut}_n$. It can be checked that an analogous 
procedure as in the case of the mapping class group gives rise
to certain homomorphisms
$$
\tau_k: \mathrm{Aut}\,F_n(k)\lra \mathrm{Hom}(H_{n},\mathcal{L}_{n}(k+1))
$$
and the totality $\{\tau_k\}_k$ of these homomorphims
induces an injective homomorphism
of graded Lie algebras
\begin{equation}
\mathrm{Gr}^+(\mathrm{Aut}\,F_n)=\bigoplus_{k=1}^\infty
\mathrm{Aut}\,F_n(k)/\mathrm{Aut}\,F_n(k+1)\lra \mathrm{Der}^+(\mathcal{L}_n).
\label{eq:johnson-aut}
\end{equation}
We refer to \cite{Kawazumi05a}\cite{Pettet}\cite{Sato05b} for
some of the recent works related to the above homomorphism.

\par
\vspace{1cm}

\section{A theorem of Kontsevich}

In this section, we recall a theorem of Kontsevich described in
\cite{Kontsevich93}\cite{Kontsevich94} which is the key
result for the argument given in the next section.
See also the paper \cite{CoV03} by Conant and Vogtmann
for a detailed proof as well as discussion of
this theorem in the context of cyclic operads.
In the above cited papers, Kontsevich considered
three kinds of infinite dimensional Lie algebras denoted by
$\mathfrak{c}_g, \mathfrak{a}_g, \mathfrak{l}_g$
(commutative, associative, and lie version, respectively).
The latter two Lie algebras are defined by
\begin{align*}
\mathfrak{a}_g &=\{\text{derivation $D$ of the tensor algebra
$T^*(H_\bQ)$}\\
&\hspace{6cm} \text{such that $D(\om_0)=0$} \}\\
\mathfrak{l}_g &=\{\text{derivation $D$ of the free Lie algebra
$\mathcal{L}^\bQ_{2g}\subset T^*(H_\bQ)$}\\
& \hspace{6cm} \text{such that $D(\om_0)=0$}\}.
\end{align*}
There is a natural injective Lie algebra homomorphism
$\mathfrak{l}_g \ra \mathfrak{a}_g$.
The degree $0$ part of both of $\mathfrak{a}_g, \mathfrak{l}_g$
is the Lie algebra $\mathfrak{sp}(2g,\bQ)$ of $\mathrm{Sp}(2g,\bQ)$.
Let $\mathfrak{a}_g^+$ (resp. $\mathfrak{l}_g^+$) denote 
the Lie subalgebra of $\mathfrak{a}_g$ (resp. $\mathfrak{l}_g$)
consisting of derivations with {\it positive} degrees.
Then the latter one $\mathfrak{l}_g^+$ is nothing other than the Lie 
algebra $\mathfrak{h}_{g,1}^\bQ$ considered in $\S 5$.
Now Kontsevich described the {\it stable} homology groups of 
the above three
Lie algebras (where $g$ tends to $\infty$)
in terms of cohomology groups of
graph complexes, moduli spaces $\mathbf{M}^m_g$
of Riemann surfaces
and the outer automorphism groups $\mathrm{Out}\, F_n$
of free groups
(or the moduli space of graphs), respectively. 
Here is the statement for the cases of $\mathfrak{a}_\infty,\mathfrak{l}_\infty$.

\begin{theorem}[Kontsevich]
There are isomorphisms
\begin{align*}
PH_*(\mathfrak{a}_\infty)&\cong
PH_*(\mathfrak{sp}(\infty,\bQ))\oplus \bigoplus_{g\geq 0,m\geq 1,
2g-2+m>0} H^*(\mathbf{M}_g^m;\bQ)^{\mathfrak{S}_m},\\
PH_*(\mathfrak{l}_\infty)&\cong\quad 
PH_*(\mathfrak{sp}(\infty,\bQ))\oplus\bigoplus_{n\geq 2} 
H^*(\mathrm{Out}\, F_n;\bQ).
\end{align*}
\end{theorem}

Here $P$ denotes the {\it primitive parts} of 
$H_*(\mathfrak{a}_\infty),H_*(\mathfrak{l}_\infty)$ which have
natural structures of Hopf algebras and
$\mathbf{M}_g^m$ denotes the moduli space of genus
$g$ smooth curves with $m$ punctures.

Here is a very short outline of the proof of the above theorem.
Using natural cell structure of the Riemann moduli space
${\bf M}_g^m\ (m\geq 1)$ and the moduli space ${\bf G}_n$
of graphs, which serves as the (rational) classifying space of 
$\mathrm{Out}\,F_n$ by \cite{CV},
Kontsevich introduced a natural filtration
on the cellular cochain complex of these moduli spaces.
Then he proved that the associated spectral sequence
degenerates at the $E_2$-term and only the diagonal
terms remain to be non-trivial. On the other hand,
by making use of classical representation theory for the
group $\mathrm{Sp}(2g,\bQ)$, he
constructed a quasi isomorphism between the $E_1$-term and
the chain complexes of the relevant Lie algebras
($\mathfrak{a}_\infty$ or $\mathfrak{l}_\infty$).
For details, see the original papers cited above
as well as \cite{CoV03}.

There is also the dual version of the above theorem
which connects the primitive cohomology of the relevant
Lie algebras with the homology groups of the moduli space
or $\mathrm{Out}\,F_n$. We would like to describe it in a
detailed form because this version is most suitable for our
purpose. The Lie algebras $\mathfrak{l}^+_\infty, \mathfrak{a}^+_\infty$
are graded. Hence their $Sp$-invariant cohomology groups are
{\it bigraded}. Let $H^k(\mathfrak{l}_\infty^+)^{Sp}_{n}$ 
and $H^k(\mathfrak{a}_\infty^+)^{Sp}_{n}$ denote
the weight $n$ part of $H^k(\mathfrak{l}_\infty^+)^{Sp}$
and $H^k(\mathfrak{a}_\infty^+)^{Sp}$ respectively.
Then we have the following isomorphisms.

\begin{align}
PH^k(\mathfrak{l}_\infty^+)^{Sp}_{2n}&\cong H_{2n-k}(\mathrm{Out}\, F_{n+1};\bQ)
\label{eq:lout}\\
PH^k(\mathfrak{a}_\infty^+)^{Sp}_{2n}&\cong 
\bigoplus_{2g-2+m=n} H_{2n-k}(\mathbf{M}_{g}^m;\bQ)_{\mathfrak{S}_m}
\end{align}

\par
\vspace{1cm}

\section{Constructing homology classes of $\mathrm{Out}\, F_n$}

In this section, we combine our construction of many cohomology classes
in $H^*(\mathfrak{h}^\bQ_{g,1})^{Sp}$ given in $\S 6, \S 7$ 
with Kontsevich's theorem given in $\S 9$ to produce
homology classes of the group $\mathrm{Out}\, F_n$.

First we see that the homomorphism given in \eqref{eq:cgh}
is {\it stable} with respect to $g$. More precisely,
the following diagram is commutative

\begin{equation}
\begin{CD}
\bQ[\Ga; \Ga\in \mathcal{G}^{odd}] @>{\Phi}>> H^*(\mathfrak{h}_{g+1,1}^\bQ)^{Sp}\\
@| @VVV\\
\bQ[\Ga; \Ga\in \mathcal{G}^{odd}] @>>{\Phi}> H^*(\mathfrak{h}_{g,1}^\bQ)^{Sp},
\end{CD}
\label{eq:phiphi}
\end{equation}
where the right vertical map is induced by the inclusion
$\mathfrak{h}_{g,1}\ra\mathfrak{h}_{g+1,1}$.
This follows from the fact that the traces $\mathrm{trace}(2k+1)$
as well as $\tau_1$ are all stable with respect to $g$ in an
obvious way. 

Next, we see that the cohomology class 
$\Phi(\Ga) \in H^*(\mathfrak{h}_{\infty,1}^\bQ)^{Sp}$
obtained in this way is primitive if and only if $\Ga$ is connected.
Keeping in mind the fact $\mathfrak{h}_{\infty,1}^\bQ=\mathfrak{l}_\infty^+$,
the property \eqref{eq:phiga},
as well as the version of Kontsevich's theorem given in
\eqref{eq:lout}
we now obtain the following theorem.

\begin{theorem}
Associated to each connected graph $\Ga\in \mathcal{G}^{odd}$ 
whose valencies are $v^a_3$ times $3$ (alternating), $v^s_3$ times $3$
(symmetric) and $v_{2k+1}$ times $2k+1\ (2k+1>3)$, we have 
a homology class
$$
\Phi(\Ga)\in H_{2n-d}(\mathrm{Out}\,F_{n+1};\bQ)
$$
where 
$$
d=v^a_3+v^s_3+v_5+v_7+\cdots,\
2n=v^a_3+3 v^s_3+5 v_5+7 v_7+\cdots.
$$
\end{theorem}

\begin{remark}
Let $\Ga_{2k+1}$ be the connected graph with two vertices
both of which have valency $2k+1$. Then $d=2, 2n=4k+2$ and
$\Phi(\Ga_{2k+1})\in H_{4k}(\mathrm{Out}\,F_{2k+2};\bQ)$
is the class $\mu_{k}$ already mentioned in $\S 4$.
\end{remark}

The following problem is an enhancement of Problem \ref{pb:grphi}
in the context of the homology of the moduli space of graphs
rather than the cohomology of the Lie algebra $\mathfrak{h}_{g,1}^\bQ$.

\begin{problem}
Give examples of odd valent graphs $\Ga$ whose associated
homology classes $\Phi(\Ga)\in H_*(\mathrm{Out}\,F_n;\bQ)$ 
are non-trivial as many as
possible. Also compare these classes with the homology classes
constructed by Conant and Vogtmann \cite{CoV04} as explicit cycles in the moduli
space of graphs. Furthermore investigate whether these classes survive
in $H_*(\mathrm{GL}(n,\bZ);\bQ)$, or else come from
$H_*(\mathrm{IOut}_n;\bQ)$, or not.
\end{problem}

\begin{remark}
It seems that the geometric meaning of Kontsevich's theorem
is not very well understood yet. In particular, there is almost
no known relations between the classes $\Phi(\Ga)\in H_*(\mathrm{Out}\,F_n;\bQ)$ 
where the rank $n$ varies. However, there should be some unknown
structures here. For example, there are graphs $\Ga$ whose associated
classes $\Phi(\Ga)$ lie in $H_8(\mathrm{Out}\,F_n;\bQ)$ for $n=7,8$
which might be closely related to the class 
$\mu_2\in H_8(\mathrm{Out}\,F_6;\bQ)$.
\end{remark}

\par
\vspace{1cm}

\section{Group of homology cobordism classes of homology cylinders}

In his theory developed in \cite{Habiro},
Habiro introduced the concept of {\it homology cobordism of surfaces}
and proposed interesting problems concerning it. 
Goussarov \cite{Goussarov} also studied the same thing in his theory.
It played an important role in the classification theory of $3$-manifolds
now called the  Goussarov-Habiro theory.
Later Garoufalidis and Levine \cite{GL} and Levine \cite{Levine}
used this concept to define a group $\mathcal{H}_{g,1}$
which consists of homology cobordism classes of {\it homology cylinders}
on $\Sg\setminus\mathrm{Int}\, D$.
We refer to the above papers for the definition
(we use the terminology {\it homology cylinder} following them)
as well as many interesting questions concerning the structure
of $\mathcal{H}_{g,1}$. It seems that the importance of this group
is growing recently.

Here we summarize some of the results of \cite{GL}\cite{Levine}
which will be necessary for our purpose here. 
Consider $\vGa=\pi_1 (\Sg\setminus\mathrm{Int}\, D)$, which is isomorphic
to $F_{2g}$, and let $\{\vGa_k\}_k$ be its lower central series as before.
Note that $\vGa$ contains a particular element $\ga\in \vGa$ which corresponds
to the unique relation in $\pi_1\Sg$.
They define
\begin{align*}
\mathrm{Aut}&_0 (\vGa/\vGa_k)=\{f\in \mathrm{Aut}(\vGa/\vGa_k);\\
&f\ \text{lifts to an endomorphism of}\ \vGa\
\text{which fixes}\ \ga\ \mathrm{mod}\ \vGa_{k+1}\}.
\end{align*}
By making a crucial use of a theorem of Stallings \cite{Stallings},
for each $k$ they obtain a homomorphism
$$
\sigma_k: \mathcal{H}_{g,1}\lra \mathrm{Aut}_0 (\vGa/\vGa_k).
$$
The following theorem given in \cite{GL} is a basic result 
for the study of the structure of the
group $\mathcal{H}_{g,1}$.

\begin{theorem}[Garoufalidis-Levine]
The above homomorphism $\sigma_k$ is surjective for any $k$.
\end{theorem}

They use the homomorphisms $\{\sigma_k\}_k$ to define a certain
filtration $\{\mathcal{H}_{g,1}(k)\}_k$ of $\mathcal{H}_{g,1}$
and show that the Johnson homomorphisms are defined also
on this group. Furthermore they are {\it surjective} so that
there is an isomorphism  
$$
\mathrm{Gr}^+(\mathcal{H}_{g,1})=
\bigoplus_{k=1}^\infty
\mathcal{H}_{g,1}(k)/\mathcal{H}_{g,1}(k+1)\cong \mathfrak{h}_{g,1}.
$$
They concluded from this fact that $\mathcal{H}_{g,1}$ contains
$\M1g$ as a subgroup. 
The homomorphisms $\sigma_k$ fit together to define
a homomorphism
$$
\sigma: \mathcal{H}_{g,1}\lra \varprojlim \mathrm{Aut}_0(\vGa/\vGa_k).
$$
They show that this homomorphism is not surjective
by using the argument of Levine \cite{Levine89}.
We mention a recent paper \cite{Sakasai05} of Sakasai for a related work.
Also they point out that,
although the restriction of $\sigma$ to $\M1g$ is injective,
$\sigma$ has a rather big kernel
because $\mathrm{Ker}\,\sigma$ at least contains the group
$$
\Theta^3_\bZ=\{\text{oriented homology $3$-sphere}\}
/\text{homology cobordism}.
$$
It is easy to see that $\Theta^3_\bZ$ is contained in the 
center of $\mathcal{H}_{g,1}$ so that we have a central extension
\begin{equation}
0\lra \Theta^3_\bZ\lra \mathcal{H}_{g,1}\lra \overline{\mathcal{H}}_{g,1}\lra 1
\label{eq:hhc}
\end{equation}
where $\overline{\mathcal{H}}_{g,1}=\mathcal{H}_{g,1}/\Theta^3_\bZ$.

The group $\Theta^3_\bZ$ is a
very important abelian group in low dimensional topology. 
In \cite{Furuta},
Furuta first proved that this group is 
an infinitely generated group by making use of gauge theory. 
See also the paper \cite{FS} by Fintushel and Stern for another proof.
However, only a few additive invariants are known on
this group at present besides 
the classical surjective homomorphism
\begin{equation}
\mu: \Theta^3_\bZ\lra \bZ/2
\label{eq:Rokhlin}
\end{equation}
defined by the Rokhlin invariant.  One is 
a non-trivial homomorphism $\Theta^3_\bZ\lra \bZ$
constructed by  Fr$\phi$yshov \cite{Froyshov} and
the other is given by Ozsv\'ath and Szab\'o \cite{OZ}
as an application of their Heegaard Floer homology theory.
Neumann \cite{Neumann} and Siebenmann \cite{Siebenmann} defined
an invariant $\bar\mu$ for plumbed homology $3$-spheres and 
Saveliev \cite{Saveliev99} introduced his $\nu$-invariant for
any homology sphere by making use of the Floer homology.
On the other hand, Fukumoto and Furuta (see \cite{FF}\cite{FFU}) defined
certain invariants for plumbed homology $3$-spheres using gauge theory.
According to a recent result of Saveliev \cite{Saveliev02},
these invariants fit together to give a candidate of another homomorphism
on $\Theta^3_\bZ$.

The situation being like this, $\Theta^3_\bZ$ 
remains to be a rather mysterious group.
Thus we have the following important problem.

\begin{problem}
Study the central extension \eqref{eq:hhc} from the point of
view of group cohomology as well as geometric topology. In particular
determine the Euler class of this central extension which is 
an element of the group
$$
H^2(\overline{\mathcal{H}}_{g,1};\Theta^3_\bZ)\cong
\mathrm{Hom}(H_2(\overline{\mathcal{H}}_{g,1}),\Theta^3_\bZ)
\oplus \mathrm{Ext}_\bZ(H_1(\overline{\mathcal{H}}_{g,1}),\Theta^3_\bZ).
$$
\end{problem}

This should be an extremely difficult problem.
Here we would like to indicate a possible method of attacking it,
and in particular a possible way of
obtaining additive invariants for the group 
$\Theta^3_\bZ$, very briefly. Details will be given in a 
forthcoming paper.

Using the traces, we can define a series of 
certain cohomology classes
$$
\tilde{t}_{2k+1}\in H^2(\varprojlim \mathrm{Aut}_0(\vGa/\vGa_m))\quad (k=1,2,\cdots).
$$
These are the ``group version" of the elements 
$t_{2k+1}\in H^2(\mathfrak{h}_{g,1}^\bQ)^{Sp}$
defined in $\S 6$ (Definition \ref{def:et}).
The homomorphism $\sigma$ is trivial on $\Theta^3_\bZ$ so that
we have the induced homomorphism
$\bar\sigma: \overline{\mathcal{H}}_{g,1}\ra 
\varprojlim \mathrm{Aut}_0(\vGa/\vGa_m)$.

\begin{conjecture}
\begin{align*}
&\text{
1. $\bar\sigma^*(\tilde{t}_{2k+1})$ is non-trivial in 
$H^2(\overline{\mathcal{H}}_{g,1})$
for any $k$}\\
&\text{
2. $\sigma^*(\tilde{t}_{2k+1})$ is trivial in
$H^2(\mathcal{H}_{g,1})$ for any $k$}.
\end{align*}
\end{conjecture}

The first part of the above conjecture is the 
``group version" of Conjecture \ref{conj:t}
and it should be even more difficult to prove.
On the other hand, if the classes $\sigma^*(\tilde{t}_{2k+1})$
were non-trivial, then they would serve as invariants for
certain $4$-manifolds ($2$-dimensional family of homology
cylinders). This seems unlikely to be the case.
Thus the second part is related to the following problem.

\begin{problem}
Determine the abelianization of the group $\mathcal{H}_{g,1}$.
Is it trivial? Also determine the second homology group
$H_2(\mathcal{H}_{g,1};\bZ)$. Is the rank of it equal to
$1$ given by the signature?
\end{problem}

If everything will be as expected,
we would obtain non-trivial homomorphisms
$$
\hat{t}_{2k+1}: \Theta^3_\bZ\lra \bZ
$$
as {\it secondary} invariants
associated to the cohomology classes $\tilde t_{2k+1}$.
There should be both similarity and difference between these cases
and the situation where we interpreted
the Casson invariant as the secondary invariant associated
to the first Mumford-Morita-Miller class $e_1$ (see \cite{Morita89}).
More precisely, they are similar because they are all related to
some cohomology classes in $H^2(\mathfrak{h}^\bQ_{g,1})$.
They are different because $e_1$ is non-trivial in
$H^2(\mathcal{H}_{g,1})$ whereas we expect that the other classes
$\sigma^*(\tilde{t}_{2k+1})$ would be all trivial in the same group.

Also recall here that Matumoto \cite{Matumoto} and Galewsky and Stern \cite{GS}
proved that every topological manifold (of dimension $n\geq 7$)
is simplicially triangulable if and only if 
the homomorphism \eqref{eq:Rokhlin} splits.
In view of this result, it should be important to investigate
the mod $2$ structure of the extension \eqref{eq:hhc}
keeping in mind the works of Birman and Craggs \cite{BC} as well as
Johnson \cite{Johnson85b}
for the case of the mapping class group.

However we have come too far and surely many things have to be
clarified before we would understand the structure of the group
$\mathcal{H}_{g,1}$.

Finally we would like to propose a problem.

\begin{problem}
Generalize the infinitesimal presentation of the Torelli Lie algebra given 
by Hain \cite{Hain97} to the case of the group of homology cobordism
classes of homology cylinders.
\end{problem}

\par
\vspace{1cm}

\section{Diffeomorphism groups of surfaces}

Let us recall the important problem of
constructing (and then detecting) 
characteristic classes of smooth fiber bundles as well
as those of foliated fiber bundles whose fibers are
diffeomorphic to a general closed $C^\infty$ manifold $M$.
This is equivalent to the problem of computing the cohomology groups
$H^*(\mathrm{BDiff}\, M)$ (resp. $H^*(\mathrm{BDiff}^\delta M)$)
of the classifying space of the diffeomorphim group $\mathrm{Diff}\,M$
(resp. the same group $\mathrm{Diff}^\delta M$ equipped with
the discrete topology) of $M$.
Although the theory of higher torsion invariants of fiber bundles 
developed by Igusa \cite{Igusa} on the one hand and by Bismut, Lott \cite{BL}
and Goette on the other and also the Gel'fand-Fuks cohomology $H^*_{GF}(M)$ of $M$
(see e.g. \cite{Haefliger})
produce characteristic classes for the above two types of $M$-bundles,
there seem to be only a few known results concerning explicit computations
for specific manifolds. The same problems are also important for various
subgroups of $\mathrm{Diff}\,M$, in particular the symplectomorphism
group $\mathrm{Symp}(M,\om)$ (in the case where there is given 
a symplectic form $\om$ on $M$)
as well as the volume preserving diffeomorphism group.

Here we would like to propose two problems
for the special, but
at the same time very important case of surfaces.
Note that we have two characteristic classes
\begin{equation}
\int_{fiber}u_1c_1^2, \int_{fiber}u_1c_2\in H^3({\rm BDiff}_+^\delta\Sg;\bR)
\label{eq:uc}
\end{equation}
which are defined to be the fiber integral of the characteristic classes
$u_1c_1^2, u_1c_2$
of codimention two foliations.
In the case of $g=0$ (namely $S^2$),
Thurston and then Rasmussen \cite{Rasmussen}
(see also Boullay \cite{Boullay}) proved that these classes are linearly
independent and vary continuously.
Although it is highly likely that their results can be extended to the
cases of surfaces of higher genera, explicit construction seems to be
open. 

\begin{problem}
Prove that the above characteristic classes induce surjective 
homomorphism
$$
H_3(\mathrm{BDiff}^\delta_+\Sg;\bZ)\lra \bR^2
$$
for any $g$.
\end{problem}

The cohomology classes in \eqref{eq:uc} are stable with respect to $g$.
On the other hand, in \cite{KM05a}\cite{KM05} we found an
interesting interaction between the
twisted cohomology group of the mapping class group and 
some well known concepts in symplectic topology
such as the flux homomorphism as well as the Calabi homomorphism
 (see \cite{MS} for generalities of the symplectic topology). 
By making use of this, we defined certain
cohomology classes of $\mathrm{BSymp}^\delta\Sg$ and proved
non-triviality of them.

In view of the fact that all the known cohomology classes of
$\mathrm{BDiff}^\delta_+\Sg$ as well as $\mathrm{BSymp}^\delta\Sg$
are stable with respect to the genus, we would like to
ask the following problem.

\begin{problem}
Study whether the homology groups of $\mathrm{BDiff}^\delta_+\Sg$
stabilize with respect to $g$ or not.
The same problem for the group $\mathrm{Symp}^\delta\Sg$.
\end{problem}
\par
\vspace{5mm}
\noindent
{\it Acknowledgments}\quad
The author would like to express his hearty
thanks to R. Hain, N. Kawazumi, D. Kotschick, M. Matsumoto, 
H. Nakamura, T. Sakasai for enlightening discussions as well as useful informations
concerning the problems treated in this paper.

\par
\vspace{1cm}

\bibliographystyle{amsplain}

\end{document}